\documentclass{article}
\usepackage[latin1]{inputenc}
\usepackage{t1enc}
\usepackage{a4}
\usepackage{amsmath}
\usepackage{amssymb}
\usepackage{amsmath,theorem,epsf,psfrag}
\usepackage[pdftex]{graphicx}
\newtheorem{theorem}{Theorem}
\newtheorem{definition}{Definition}
\newtheorem{lemma}{Lemma}
\newtheorem{corollary}{Corollary}

\newtheorem{remark}{Remark}

\newcommand{\EE}{\mathbb{E}}
\newcommand{\ZZ}{\mathbb{Z}}

\newcommand{\NN}{\mathbb{N}}
\newcommand{\PP}{\mathbb{P}}
\newcommand{\QQ}{\mathbb{Q}}
\newcommand{\RR}{\mathbb{R}}

\newcommand{\TT}{\mathbb{T}}

\newcommand{\hH}{{\cal H}}
\newcommand{\hF}{{\cal F}}
\newcommand{\hB}{{\cal B}}
\newcommand{\hP}{{{\cal P}_d}}

\newcommand{\eps}{{\varepsilon}}

\begin{document}
\begin{center}
{\Large
 Tessellation-valued processes that are generated by cell division} 

\vspace{0.5cm}

{Werner Nagel}\\
{Friedrich-Schiller-Universit\"at Jena,\\
Institut f\"ur Stochastik,\\
Ernst-Abbe-Platz 2,
07743 Jena, Germany.\\
Email: werner.nagel@uni-jena.de}  

\vspace{0.5cm}

{Servet Mart{\'i}nez}\\
{Universidad de Chile,
Departamento Ingenier{\'i}a Matem\'atica and Centro
Modelamiento Matem\'atico, \\
UMI 2807 CNRS, Casilla 170-3, Correo 3,\\
Santiago, Chile.  \\
Email: smartine@dim.uchile.cl}

\end{center}

\begin{abstract}
\noindent Processes of random tessellations of the Euclidean space $\mathbb{R}^d$, $d\geq 1$, are considered which are generated by subsequent division of their cells. Such processes are characterized by the laws of the life times of the cells until their division and by the laws for the random hyperplanes that divide the cells at the end of their life times. 
The STIT tessellation processes are a reference model. In the present paper a generalization concerning the life time distributions is introduced, a sufficient condition for the existence of such cell division tessellation processes is provided and a construction is described. In particular, for the case that the random dividing hyperplanes have a Mondrian distribution -- which means that all cells of the tessellations are cuboids -- it is shown that the intrinsic volumes, except the Euler characteristic, can be used as the parameter for the exponential life time distribution of the cells.
\end{abstract}

\vspace{.5cm}


\vspace{.5cm}


\section{Introduction}

The study of random tessellations (or mosaics) is a substantial part of stochastic geometry. The Voronoi tessellations and Poisson hyperplane tessellations are classical and well-established. Motivated by modeling fracture or crack patterns as they are observed for example in geology, materials science, nanotechnology or drying soil patterns \cite{xia_hutchinson_2000, Boulogne_et_al_2015, Goehring_et_al2016, Hafver_et_al2014,  seghir/arscott15}, also space-time processes of tessellations are considered where the cells are consecutively divided. 
Recently, such models are also of growing interest in the context of machine learning, see \cite{OReilly_Tran2022} and the references therein.

Several of the approaches suggested so far, turned out to be cumbersome, if not unfeasible, with respect to theoretical investigations, and only simulation studies can be performed.
In \cite{nagel_weiss2005} the STIT tessellation process was introduced, which is a cell division process allowing for numerous theoretical results, see for example \cite{nagel_weiss2008, mecke_nagel_weiss2007, mecke_nag_weiss2008_global_construction, schreiber_thae2012_second_order, schreiber_thae2013_limit_theorems, schreiber_thae2013_connection_PHT, martinez_nagel20016}. Due to its nice mathematical properties, it can be considered as a reference model for fracture patterns. On the other hand, some statistical goodness-of-fit checks indicate the need for modification and adaption of the STIT model, see \cite{nagel_mecke_ohser_weiss2008, leon_ohser_nagel_arscott2019}.
This motivates the present paper.

Our work is essentially inspired by Cowan's paper \cite{cowan2010} where a systematic approach to a wide class of cell division processes is provided by introducing the 'selection rules' -- in the present paper adapted as 'life time distributions' -- and 'division rules'.

In \cite{schreiber_thae2013_shape_driven, georgii_schreiber_thaele2015} a theoretical base for generalizations of the STIT model was provided. In those papers the focus was on modifications of the division rules for the cells. 

Our purpose is the construction of tessellation-valued cell division processes in the $d$-dimensional Euclidean space $\RR^d$, where the life time distributions of the cells differ from those ones of the STIT model. 
We provide sufficient conditions and a rigorous proof for their existence, which is a critical issue for such models, and we study some of their features.

The idea of the proof of Theorem \ref{thm: sufficient existence} is based on the 'global construction' of STIT tessellations by Mecke in \cite{mecke_nag_weiss2008_global_construction}. For this reason, the models which we consider, all have the same division rule like the STIT tessellation, driven by a translation invariant hyperplane measure. The new results of this paper concern a variety of life time distributions of the cells.

In Section \ref{sec:zero-cells} we start with a proof of the tail triviality of the zero-cell process of a Poisson hyperplane tessellation process. Because this process of zero-cells plays a central role in several proofs of cell division processes, this is a result of its own interest.
Then in Section \ref{sec:cell division process} a sufficient condition for the existence of a class of cell division processes is provided together with their construction.
Several more concrete results are shown in Section \ref{sec:Mondrian} for the particular case of the so-called Mondrian model where all the tessellation cells are cuboids.
We conclude by describing a link between a particular cell division process in a bounded window and a fragmentation as defined in \cite{bertoin2006}.

\section{Notation}\label{section:Notation}

By	$\RR^d$  we denote the $d$-dimensional Euclidean space, $d\geq 1$, with the scalar product $\langle \cdot ,\cdot \rangle$ .  The origin	is $o$, and $S^{d-1}$ is the unit sphere with the Borel $\sigma$-algebra ${\mathcal B} (S^{d-1})$. Denote by ${\sf B}_r$ the ball centered in the origin $o$ with radius $r>0$. The topological interior of a set $A\subseteq \RR^d$ is ${\rm int} (A)$. The Lebesgue measure on $\RR$ is denoted by $\lambda$.
The set of the nonnegative integers is $\NN_0$, $\ZZ$ denotes the set of integers and $\ZZ_- :=\{ i\in \ZZ :\, i\leq 0\}$.
The symbol $\stackrel{D}{=}$ means the identity of distributions of two random variables. We use the abbreviation 'i.i.d.' for 'independent identically distributed' random variables.

Let $\hH $ denote the set of all hyperplanes in $\RR^d$, $\hH $ endowed with the Borel $\sigma$-algebra associated to the Fell topology.
For a hyperplane  $h\in \hH$ we use the parametrization 
$$h(u,x):= \{ y\in \RR^d:\, \langle y,u\rangle =x\}, \quad u\in S^{d-1},\, x\in \RR . $$
Note, that $h(u,x)=h(-u,-x)$, thus all hyperplanes can be parametrized in two different ways.

Let $\varphi$  denote an even probability measure on $S^{d-1}$, which is a probability measure on $S^{d-1}$ with $\varphi (A)=\varphi (-A)$ for all $A\in {\mathcal B} (S^{d-1})$.
 
Throughout this paper we use a translation invariant measure $\Lambda$ on $\hH$ which satisfies   
\begin{equation}\label{eq:hyperplane measure}
\int_{\hH} f  {\rm d} \Lambda  = \int_{S^{d-1}} \int_\RR f(h(u,r))\lambda ({\rm d} r)\varphi ({\rm d} u)
\end{equation}
for all nonnegative measurable functions  $f:\hH \to \RR$.
The probability measure $\varphi$ is called the spherical directional distribution.
We need a further assumption on $\varphi$ which guarantees that the constructed structures are indeed tessellations with bounded cells, see p. 486 of \cite{sw08}, 
\begin{equation}\label{assumption:hyperplane measure}
\mbox{$\varphi$ is not concentrated on a great subsphere of } S^{d-1}.
\end{equation}
A great subsphere is the intersection of $S^{d-1}$ with a hyperplane through the origin. Condition \eqref{assumption:hyperplane measure} is equivalent to the assumption that there is no line in $\RR^d$ to which $\varphi$-almost-all hyperplanes are parallel.

A polytope in $\RR^d$ is the convex hull of a nonempty finite set of points. By $\hP$  we denote the set of all $d$-dimensional polytopes in $\RR^d$.

\begin{definition}\label{def:tessellation}
A {\em tessellation of $\RR^d$} is a countable set $T\subset \hP$ of $d$-dimensional polytopes satisfying the following conditions:
\begin{enumerate}
\item[(a)] $ \displaystyle{\bigcup_{z\in T} z = \RR^d }$ (covering),
\item[(b)] If $z,z'\in T$ and $z\not= z'$, then ${\rm int}(z)\cap {\rm int}(z')= \emptyset$ (disjoint interiors),
\item[(c)] The set $\{ z\in T: z\cap C \not= \emptyset \}$ is finite for all compact sets $C\subset \RR^d$ (local finiteness).
\end{enumerate}
\end{definition}

The polytopes $z\in T$ are referred to as the {\em cells} of $T$.

Let $\TT$ denote the set of all tessellations of $\RR^d$, endowed with the usual $\sigma$-algebra, which is associated with the Borel $\sigma$-algebra for the Fell topology in the space of closed subsets of $\RR^d$, using that the union of cell boundaries $\partial T := \bigcup_{z\in T} \partial z$, also referred to as the $(d-1)$-skeleton, of a tessellation $T$ is a closed subset 
of $\RR^d$.
A random tessellation is a measurable mapping from a probability space to $\TT$.

For more details on the definition of random tessellation we refer to \cite{sw08}.

In this paper, Poisson point processes are used extensively. For this,we refer to the book \cite{last_penrose2017}, in particular to Chapters 5 and 7 there.

\section{An auxiliary Poisson process of hyperplanes that are marked with birth times, and the process of zero-cells}\label{sec:zero-cells} 

Let $\Lambda$ be a translation invariant measure on the space $\hH $ of hyperplanes in $\RR^d$, satisfying \eqref{eq:hyperplane measure} and \eqref{assumption:hyperplane measure}. We denote by $\hat X^*$ a Poisson point process on $\hH \times (0,\infty )$ with the intensity measure $\Lambda \otimes \lambda_+$, where $\lambda_+$ is the Lebesgue measure on $(0,\infty )$. If $(h,t)\in \hat X^*$, we interpret this as a hyperplane with birth time $t$. Furthermore, define $\hat X^*_t := \{  (h,t')\in \hat X^*:\, t'\leq t\}$ and  $\hat X_t := \{ h\in \hH: \exists (h,t')\in X^*_t\}$. Now we can introduce the process $(\tilde z^o_t , t>0)$, where $\tilde z^o_t$ is the zero-cell -- that is the cell containing the origin $o$ -- of the Poisson hyperplane tessellation of $\RR^d$ induced by $\hat X_t$.
The process $(\tilde z^o_t , t>0)$ is a pure jump Markov process, with  and $\tilde z^o_{t'} \subseteq \tilde z^o_{t}$ if $t<t'$.
And, most importantly,
\begin{equation}\label{eq: space filling}
\bigcup_{t>0} \tilde z^o_t =\RR^d \quad a.s.
\end{equation}

Denote by $(t_i ,\, i \in \ZZ )$ the ordered sequence of jump times of $(\tilde z^o_t , t>0)$ with $t_0<1<t_1$. Correspondingly, $(h_i, t_i)\in \hat X^*$ are the time marked dividing hyperplanes causing a jump of $(\tilde z^o_t , t>0)$.

Denote by
\begin{equation}\label{sequ Poiss zero}
(\tilde z^o_{(i)},\, i\in \ZZ)
\end{equation}
the sequence of zero cells with 
 $\tilde z^o_{(i)}:=\tilde z^o_{t_i}$ . 

Consider the process $(\hat z_n ,\, n\in \NN_0)$ with 
$\hat z_n:= \tilde z^o_{(-n)}$,
and denote by $\hF_n := \sigma (\hat z_n)$, $n\in \NN_0$, the $\sigma$-algebra generated by $\hat z_n$\,.
The tail-$\sigma$-algebra is then 
\begin{equation*}
\hF_\infty := \bigcap_{m\in \NN_0}\bigvee_{k\geq m} \hF_k ,
\end{equation*}
where $\bigvee$ stands for the $\sigma$-algebra generated by the union of the respective $\sigma$-algebras.

\begin{lemma}\label{lemma:trivial tail zero cell}
The tail-$\sigma$-algebra  $\hF_\infty $ of the process $(\hat z_n ,\, n\in \NN_0)$ is trivial, that is $\PP(B)\in \{ 0,1\}$ for all $B\in \hF_\infty $\,.
\end{lemma}
{\bf Proof}
Let be $B\in \hF_\infty$ and $\eps >0$ fixed. Then there are an $m\in \NN$ and a $B_m \in \bigvee_{k\leq m} \hF_k $ such that $\PP (B\Delta B_m )<\eps$\,, see Theorem D, Chapter III of \cite{halmos1970}. In the following we consider $m$ and $B_m$ as fixed. 

The main idea of the proof is to construct times ${\tilde t} < t^*$, two polytopes $W'\subset W$ containing the origin $o$, and events $E$, $E_m$, $A$, $A_n$, $n>m$, such that

-- the event $E\cap E_m$ implies that $W'\subseteq \hat z_m \subseteq W$, and that the jump time $t_{-m}\geq t^*$\,,

-- the event $A\cap A_n$ implies that $W \subseteq \hat z_n$, and that the jump time $t_{-n}\leq {\tilde t}$\,.

For $W' , W\in \hP$, $W'\subset W$, with the origin $o$ in the interior of $W'$, define the time $S(W',W):= \inf \{ t>0:\, \tilde z^o_t \subseteq W \}$ which is the time of first separation of the boundary of $W'$ from the boundary of $W$ by the boundary of the zero-cell, in \cite{martinez_nagel2014} introduced as the 'encapsulation time'.

For  $t^* \in (0,1)$ and $r>0$ define the intervals $a_i:=\left((i-1)\frac{r}{m} ,\,  i\frac{r}{m}\right)$ of length $r/m$, and time intervals $\Delta_i :=\left(t^* + (i-1)\frac{1-t^*}{m} ,\, t^* + i\frac{1-t^*}{m} \right)$ of length $(1-t^*)/m$ between $t^*$ and 1,  $i\in \{ 1,\ldots ,m\}$.

Now, for $m\in \NN$ and $\eps >0$ choose a pair $(t^*,r)\in (0,1)\times (0,\infty )$ such that for the event 
\begin{align*}
E_m 
:= &\{ \hat X^* \cap 
\left(\{ h(u,x)\in \hH :\, x\in (-r,r)\}\times (0,t^*) \right)= \emptyset \} \\
& \cap \bigcap_{i=1}^m \left\{ \hat X^* \cap \left(
\{ h(u,x)\in \hH :\, x\in a_i \cup (-a_i)\} \times \Delta_i \right) \not= \emptyset \right\}
\end{align*}
holds
\begin{equation}\label{eq:cond Em}
\PP(E_m)={\rm exp}(-t^* r)\, \left( 1-{\rm exp} \left(-\frac{2(1-t^*)  \, r}{m^2}\right) \right)^m >1-\eps \, .
\end{equation}
This can be realized by choosing 
$r> -\frac{m^2}{2}\ln \left( 1-(1-\eps )^{\frac{1}{2m}} \right)$
and then $t^* < - \frac{1}{2r} \ln (1-\eps )$, which guarantees that each of the two factors is greater than $\sqrt{1-\eps}$.

The event $E_m$ means that the ball ${\sf B}_r$ is not intersected by a hyperplane of $\hat X^*$ until time $t^*$ and then in each of the $m$ time intervals\\
 $\Delta_i :=\left(t^* + (i-1)\frac{1-t^*}{m} ,\, t^* + i\frac{1-t^*}{m} \right)$ of length $(1-t^*)/m$ between $t^*$ and 1 there is at least one hyperplane of $\hat X^*$ with a distance to the origin in the interval 
$a_i:=\left((i-1)\frac{r}{m} ,\,  i\frac{r}{m}\right)$, $i\in \{ 1,\ldots ,m\}$. This guarantees that there are at least $m$ jumps of the zero-cell process in the time interval $(t^* ,1)$. 

The pair $(t^*, r)$ is now considered to be fixed.
Choose $W'\in \hP$ such that ${\sf B}_r \subset W'$.
By Lemma 5 of of \cite{martinez_nagel2014} there exists a $W\in \hP$ 
 with $W'\subset W$ such that for  the event $E:=\{ S(W',W)<t^*\}$ holds
\begin{equation}\label{eq:encaps time}
\PP(E) =\PP(S(W',W)<t^*)>1-\eps \, .
\end{equation}

From now on, also $W$ is considered to be fixed.
Furthermore, choose ${\tilde t}\in (0,t^*)$  such that for the event $A:=\{ [W] \cap \hat X_{{\tilde t}} =\emptyset  \}$, which means that at time $\tilde t$ the window $W$ is contained in the zero-cell,
\begin{equation*}
\PP(A)= {\rm exp}\left(-{\tilde t}\, \Lambda ([W])\right) >1-\eps \, .
\end{equation*}
Then, having fixed $\tilde t$, choose $n>m$ such that for the event $A_n$, that in the time interval $(\tilde t , 1)$ at most $n$ hyperplanes of $\hat X^*$ intersect $W$, 
\begin{equation*}
\PP(A_n)= \sum_{i=0}^n \frac{((1-{\tilde t}) \, \Lambda ([W]) )^i}{i!} {\rm exp}\left(-(1-{\tilde t})\, \Lambda ([W])\right) >1-\eps \, .
\end{equation*}
Note that $A\cap A_n \subseteq \{ W\subseteq \hat z_n\}$\,.

By the characteristic property of Poisson point processes, the restricted point processes on disjoint sets, namely $\hat X^* \cap ([W] \times (t^*,1))$ and
$\hat X^* \cap ([W]^c \times (0,{\tilde t}))$, are independent.
Hence, for $B_m \in \bigvee_{k\leq m} \hF_k $ and $C \in \bigvee_{k\geq n} \hF_k $ and with $D:=E\cap E_m \cap A \cap A_n$
we have that  the events $B_m$ and $C$ are conditionally independent under the condition $D$ and thus
\begin{align*}
  & \left| \PP(B_m \cap C) -\PP(B_m)\PP(C) \right|  \\
= & \left| \PP(B_m \cap C\cap D) +\PP(B_m \cap C\cap D^c)-\PP(B_m)\PP(C) \right|  \\
\leq & \left|\PP(B_m \cap D) \PP(C\cap D)/\PP(D) -\PP(B_m)\PP(C) \right|  +\PP(B_m \cap C\cap D^c) \\
< & \frac{8\, \eps + 16\, \eps^2}{1-4\eps}\, .
\end{align*}
Because $B\in \hF_\infty \subseteq \bigvee_{k\geq n} \hF_k $ 
we have $\left| \PP(B_m \cap B) -\PP(B_m)\PP(B) \right| <\frac{8\, \eps + 16\, \eps^2}{1-4\eps} \,$. Then  $\PP (B\Delta B_m )<\eps$ yields
$\left| \PP(B \cap B) -\PP(B)\PP(B) \right| <\frac{8\, \eps + 16\, \eps^2}{1-4\eps} + 2\, \eps$ for all $\eps >0$, and hence $P(B)=0$ or 1.
\hfill $\Box$

The tail triviality of the zero-cell process is now used to show further 0-1 probabilities.

\begin{corollary}\label{corollary: G tau 0-1}
Let $\Lambda$ be a translation invariant measure on the space of hyperplanes in $\RR^d$, satisfying \eqref{eq:hyperplane measure} and \eqref{assumption:hyperplane measure}, and let $G: \hP \to [0,\infty )$ be a measurable functional on the set of $d$-dimensional polytopes. 

(a) Then 
\begin{equation}\label{eq: 0 1 G}
\PP\left(\exists m\in \ZZ:\ \sum_{i\leq m} G\left(\tilde z^o_{(i)}\right)^{-1}  
                < \infty  \right) \in \{ 0,\, 1\}\, .
\end{equation}

(b)
Let $(\tau'_i,\, i\in \ZZ)$ be a sequence of i.i.d. random variables, independent of $\hat X^*$. Then
\begin{equation}\label{eq:0 1 tau}
\PP\left(\exists m\in \ZZ:\  \sum_{i\leq m} G\left(\tilde z^o_{(i)}\right)^{-1} \tau'_i 
                < \infty  \right) \in \{ 0,\, 1\}\, 
\end{equation}
and 

(c)
\begin{equation}\label{eq:explosion}
\PP\left( \exists m\in \ZZ:\ \sum_{i\leq m} G\left(\tilde z^o_{(i)}\right)^{-1} \tau'_i 
                < \infty  \right) = \PP\left( \exists m\in \ZZ:\ \sum_{i\leq m} G\left(\tilde z^o_{(i)}\right)^{-1}  
                < \infty  \right)\, .
\end{equation}

\end{corollary}
{\bf Proof}
By the measurability of $G$, the event $\left\{\exists m\in \ZZ:\ \sum_{i\leq m} G\left(\tilde z^o_{(i)}\right)^{-1} < \infty \right\}$ is an element of the tail-$\sigma$-algebra  $\hF_\infty$, and we can apply Lemma \ref{lemma:trivial tail zero cell}. 

If $(\tau'_i,\, i\in \ZZ)$ be a sequence  of i.i.d. random variables  which is independent of $\hat X^*$, then this result can be extended straightforwardly to the product-$\sigma$-algebras $\hF_n \otimes \hF_n^\tau$, $n\in \NN_0$, where 
$\hF_n := \sigma (\hat z_n)$, $n\in \NN_0$ and $\hF_n^\tau:= \sigma (\tau_n)$. This yields \eqref{eq:0 1 tau}.

To prove \eqref{eq:explosion}, first note that a.s. non-explosion
of the process
\begin{equation*}
\left( \sum_{-j\leq i\leq m} G\left(\tilde z^o_{(i)}\right)^{-1} \tau'_i ,\,     j\in \NN \right)
\end{equation*}
means that $ \PP \left(\sum_{i\leq m} G\left(\tilde z^o_{(i)}\right)^{-1} \tau'_i = \infty\right)=1$. It is well-known in the theory of birth processes, see for example Proposition 13.5 of \cite{kallenberg2021}, that this is equivalent to
$ \PP \left(\sum_{i\leq m} G\left(\tilde z^o_{(i)}\right)^{-1} = \infty\right)=1$, which means that the property of non-explosion only depends on the series of the expectations $G\left(\tilde z^o_{(i)}\right)^{-1}$ of the holding times of the respective states of the process. 
Hence, for the respective complements
\begin{equation*}
\PP \left(\sum_{i\leq m} G\left(\tilde z^o_{(i)}\right)^{-1} \tau'_i < \infty\right)=0 \ \Leftrightarrow \ \PP \left(\sum_{i\leq m} G\left(\tilde z^o_{(i)}\right)^{-1}  < \infty\right)=0\, .
\end{equation*}
Because  by \eqref{eq: 0 1 G}  and \eqref{eq:0 1 tau} both probabilities can have the values $0$ or $1$ only, we obtain
\eqref{eq:explosion}.
\hfill $\Box$

In \cite{martinez_nagel2018}  the rescaled time-stationary process of zero-cells $(a^t \tilde z^o_{a^t} , t\in \RR)$ for $a>1$ was studied in more detail and it was shown that this is a Bernoulli flow with infinite entropy.

\section{Cell division processes}\label{sec:cell division process}
We consider a class of tessellation-valued Markov processes $(T_t, t>0)$ which are characterized by the laws of the random life times of the cells and the laws for the division of the cells by a random hyperplane at the end of their life times. Furthermore,  given the state of the process, conditional  independence of the future development in different cells is assumed. 
Denote by 
$[z]:= \{ h\in \hH :\, h\cap z \not= \emptyset \}$ the set of all hyperplanes that intersect the cell $z\in \hP$.

\begin{definition}\label{def:cell div proc}
Let  be $\QQ =(\QQ_{[z]} ,\, z\in \hP )$, where  $\QQ_{[z]}$ is a probability measure on the set $[z]$.  And let $G: \hP \to [0,\infty )$ be a measurable functional defined on the set of polytopes.

A random tessellation-valued process $(T_t, t>0)$ is called an (L-$G$) (D-$\QQ$) {\em cell division process}  if
for all $t>0$ and all cells $z\in T_t$:
\begin{enumerate}
\item[(L-G)] If $G(z)>0$, then the cell $z$ has a  random life time which is exponentially distributed with parameter $G(z)$ and, given $z$, conditionally  independent of $(T_{t'}, 0<t'<t)$ and of all other cells of $T_t$. If $G(z)=0$ then the life time is infinite, which means that $z$ is never divided.
\item[(D-$\QQ$)] At the end of its life time, $z$ is divided by a random hyperplane $h_z$ with the law $\QQ_{[z]}$, and also,  given $z$, the hyperplane $h_z$ is conditionally  independent of $(T_{t'}, 0<t'<t)$ and of all other cells of $T_t$ and independent of the life time of $z$.
\end{enumerate}
\end{definition}

It is not trivial to see whether such a cell division process $(T_t, t>0)$ exists for given $G$ and $\QQ$, because there is no appropriate initial tessellation at time $t=0$.
Even if it is clear how to perform the cell division dynamics in a bounded polytope  $W$, a window, a consistent extension to the whole space $\RR^d$ is by no means straightforward.

Examples of interest for such rules can be related to translation invariant hyperplane measures $\Lambda$ which satisfy \eqref{eq:hyperplane measure} and \eqref{assumption:hyperplane measure}. 

We mean by 
\begin{enumerate}
\item[(L-$\Lambda$)] \ldots that the life time of each cell $z$ is exponentially distributed with parameter $\Lambda ([z])$, 
\item[(L-$V_d$)] \ldots that the life time of each cell $z$ is exponentially distributed with parameter $V_d(z)$, the volume of $z$,
\item[(D-$\Lambda$)] \ldots that the law of the random hyperplane dividing a cell $z$ is \\ $\QQ_{[z]} =\Lambda ([z])^{-1}  \, \Lambda (\cdot \cap [z])$, which is the probability measure on $[z]$, induced by $\Lambda$.
\end{enumerate}

The STIT tessellation process driven by $\Lambda$ is a cell division process with (L-$\Lambda$) and (D-$\Lambda$).
In \cite{nagel_biehler15}  it was shown that in the class of cell division processes defined above, only the (L-$\Lambda$) (D-$\Lambda$) model has the property of spatial consistency which is sufficient for its existence. Therefore it is of interest to show the existence of further cell division processes without requiring spatial consistency.

\subsection{A sufficient condition for the existence of a cell division process in $\RR^d$}

Let $(\hP ,\hB (\hP))$ be the measurable space of $d$-dimensional polytopes, where $\hB (\hP)$ is the Borel $\sigma$-algebra with respect to the Fell topology on the space of closed subsets of $\RR^d$, see \cite{sw08}.
A functional
$G: \hP \to [0,\infty )$
is called monotone if $G(z_1)\leq G(z_2)$ for all $z_1,z_2 \in \hP$ with $z_1\subseteq z_2$. It is called invariant under translations if $G(z)=G(z+x)$ for all $z\in \hP$ and all $x\in \RR^d$.

\begin{theorem}\label{thm: sufficient existence}
Let $\Lambda$ be a translation invariant measure on the space of hyperplanes in $\RR^d$, satisfying \eqref{eq:hyperplane measure} and \eqref{assumption:hyperplane measure}, and let $G: \hP \to [0,\infty )$ be a measurable, monotone and translation invariant functional on the set of $d$-dimensional polytopes.  If for the random sequence $\left(\tilde z^o_{(i)},\, i\in \ZZ\right)$ defined in \eqref{sequ Poiss zero} holds that a.s.
\begin{equation}\label{eq: suff cond}
 \exists m\in \ZZ:\quad  \sum_{i\leq m} G\left(\tilde z^o_{(i)}\right)^{-1}  
                < \infty  
\end{equation}
then there exists a Markov tessellation-valued process $(T_t, t>0)$, which is an (L-$G$) (D-$\Lambda$) cell division process according to Definition \ref{def:cell div proc}.
\end{theorem}

{\bf Proof}
The idea of the proof is inspired by the global construction of STIT tessellations, described in \cite{mecke_nag_weiss2008_global_construction}. The key is a backward construction for $t\downarrow 0$ of the process $(z^o_t, t>0)$ of zero-cells of $(T_t, t>0)$.

The  auxiliary process $(\hat X_t, t>0) $ of Poisson hyperplane tessellations has the division law for its zero-cells which is exactly the (D-$\Lambda$) rule. Hence 
the sizes and the shapes of the zero-cells of the Poisson hyperplane tessellations can be used for the construction of $(z^o_t, t>0)$. What has to be changed is the 'clock'. The time axis has to be transformed such that the pure jump process of zero cells receives jump times according to (L-$G$), that is, the life times of the respective cells are exponentially distributed, and the parameter is the functional value of $G$ for these cells.

We use the random sequence $(\tilde z^o_{(i)},\, i\in \ZZ)$ defined in \eqref{sequ Poiss zero} for the construction, and we assign appropriate holding times to it. 

Let $(\tau'_i ,\, i\in \ZZ )$ be a sequence of i.i.d. random variables, exponentially distributed with parameter 1, and this sequence has to be independent from  $\hat X^*$. Define $(\tau_i , \,i\in \ZZ )$ as

\begin{equation}\label{eq: def clock}
\tau_i := \sum_{j\leq i-1} G(\tilde z^o_{(j)})^{-1}\, \tau'_j, \ i\in \ZZ .
\end{equation}

Now the process $(z^o_t, t>0)$ is defined by
\begin{equation}\label{eq: def zero}
z^o_t := \tilde z^o_{(i)} \quad \mbox{ for } \tau_i \leq  t < \tau_{i+1}, \ i\in \ZZ .
\end{equation}

Thus, for the existence of the process $(z^o_t,\, t>0)$ it is sufficient that a.s. there is an $i\in \ZZ$ such that $\tau_i$ is finite. By \eqref{eq:explosion} it follows that \eqref{eq: suff cond} is equivalent to it.  This implies also that $\lim_{j\to -\infty} \tau_j = 0$, a.s., and hence $\bigcup_{t>0} z^o_t = \RR^d$, by to \eqref{eq: space filling}.

Having the existence of the process $(z^o_t,\, t>0)$ of zero-cells,  the process $(T_t,\, t>0)$ is completed as follows. Denote by $\textup{cl}$ the topological closure of a subset of $\RR^d$.
At each time $\tau_i$ , when the process $(z^o_t,\, t>0)$ has a jump from $\tilde z^o_{(i-1)}$ to  $\tilde z^o_{(i)}$, an (L-$G$) (D-$\Lambda$) cell division process according to Definition \ref{def:cell div proc} is launched in  the new separated cell 
$\textup{cl}\left(\tilde z^o_{(i-1)} \setminus \tilde z^o_{(i)}\right)$.
 
It remains to show that this actually yields a tessellation $T_t$ for all $t>0$, and that the process $(T_t,\, t>0)$ is an (L-$G$) (D-$\Lambda$) cell division process.

Condition (b) of Definition \ref{def:tessellation} is obviously satisfied. To show (a), consider an arbitrary point $x\in \RR^d$. Then \eqref{eq: space filling} implies that $x\in z^o_s$ for $0<s<t_x$, where $t_x$ is the time when $x$ and $o$ are separated by a hyperplane. Thus the construction guarantees that for all $t>0$ any point $x\in \RR^d$ belongs to a the interior or to the boundary of a cell of $T_t$.

To prove property (c) of Definition \ref{def:tessellation}, consider a compact set $C\subset \RR^d$ and a time $t_0 > 0$ such that $C\subset z^o_{t_0}$. By \eqref{eq: space filling} such a time exists.
Thus for all $t\leq t_0$, there is exactly one cell in $T_t$ that has a nonempty intersection with $C$. For $t>t_0$ the the number of cells inside $z^o_{t_0}$ of the cell division process can be described by a birth process.  By the monotonicity and translation invariance of $G$, the birth rates are dominated by $G(z^o_{t_0})\times \, number\ of\ cells\ in\ G(z^o_{t_0})$. Hence, at any time $t>0$, the number of cells intersecting $C$ is a.s. finite.

Finally, to see that the process $(T_t,\, t>0)$ is an (L-$G$)  (D-$\Lambda$) cell division process, observe that the division rule for the zero-cell process is induced by the Poisson hyperplane process $\hat X$ and thus the (D-$\Lambda$) rule is realized. 
The definition in  \eqref{eq: def zero} yields that the holding time of the state $z^o_{(i)}$ is 
$\tau_{i+1} - \tau_i = G(\tilde z^o_{(i)})^{-1}\, \tau'_i$ which shows that (L-$G$) is satisfied for the zero-cell process.
For cells which are already separated from the zero-cell, the (L-$G$) and (D-$\Lambda$) rules are satisfied by the definition of the construction.
\hfill $\Box$

In the rest of this paper, we mean by an (L-$G$)  (D-$\Lambda$) process $(T_t,\, t>0)$ always the tessellation-valued process, defined in the proof.
Note that the proof does not imply that the distribution of  (L-$G$)  (D-$\Lambda$) process is uniquely determined by $G$ and $\Lambda$.

An example where the sufficient condition \eqref{eq: suff cond} is obviously not satisfied, is $G(z)=c$ for all $z\in \hP$ and some fixed $c>0$. It is not known whether a cell division tessellation process exists for this $G$, but the simulations in \cite{cowan2010} for the corresponding 'equally likely' case suggest that the cell division does not yield a tessellation according to Definition \ref{def:tessellation}, because the local finiteness condition (c) seems to be violated.

\subsection{Stationarity in space}

Now it will be shown that the translation invariance of the hyperplane measure $\Lambda$ and of the functional $G$ are sufficient for the spatial stationarity of the constructed tessellations at any fixed time $t>0$.

\begin{theorem}
If the assumptions of Theorem \ref{thm: sufficient existence} are satisfied, then for all $t>0$ the random tessellation $T_t$, described in the proof of the theorem, is a spatially stationary (or homogeneous) tessellation.
\end{theorem}
{\bf Proof}
The proof can be sketched as follows.  Because the hyperplane measure $\Lambda$ defined in \eqref{eq:hyperplane measure} is invariant under translations of $\RR^d$, the distribution of  the auxiliary Poisson process $\hat X^*$ is invariant under translations by $(x,0)\in \RR^d \times \{ 0\}$.  Thus, denoting by $\left(\tilde z^o_{(i)}(\hat X^* )  ,\, i\in \ZZ\right)$ the process of zero-cells generated by $\hat X^*$, 
\begin{equation*}
\left(\tilde z^o_{(i)}(\hat X^* +(x,0)),\, i\in \ZZ\right) \stackrel{D}{=}  \left(\tilde z^o_{(i)}(\hat X^* )  ,\, i\in \ZZ\right) \quad \mbox{ for all } x\in \RR^d .
\end{equation*}
This means the invariance under translations of $\RR^d$ of the distribution of the process
$(\tilde z^o_{(i)},\, i\in \ZZ)$  defined in \eqref{sequ Poiss zero}.

Furthermore, the translation invariance of $G$ guarantees that the life time distributions of the cells remain translation invariant too. Finally, because $\Lambda$ is invariant under translations, it follows that
$\Lambda ([z+x])^{-1}  \, \Lambda (A+x \cap [z+x]) = \Lambda ([z])^{-1}  \, \Lambda (A \cap [z])$ for all $z\in \hP$, $x\in \RR^d$ and
Borel sets $A\subset [z]$. Hence, the division law is translation-equivariant. 
\hfill $\Box$

\subsection{Construction of the cut-out appearing in a bounded window $W$}
For simulations of the model it is of interest to construct the restriction $(T_t \land W, t>0)$ of $(T_t,\, t>0)$ to a bounded window $W\in \hP$,
where $T_t \land W := \{ z\cap W : z\in T_t, \textup{int}(z\cap W)\not= \emptyset \}$. In \cite{nagel_weiss2005} it was shown, that for STIT tessellations the (L-$\Lambda$)  (D-$\Lambda$) cell division process can be launched in $W$ without regarding boundary effects, because the STIT model is spatially consistent. And in \cite{nagel_biehler15} it was shown that all the other models considered in the present paper are not spatially consistent which means that for the construction of $(T_t \land W, t>0)$  also information  outside of $W$ is needed.

The construction described in the proof of Theorem \ref{thm: sufficient existence} also shows how $(T_t \land W, t>0)$ can be realized. For a given window $W$ construct a zero-cell $\tilde z^o_{(i_W)}$ from the process $(\tilde z^o_{(i)},\, i\in \ZZ)$ defined in \eqref{sequ Poiss zero} such that $ W \subset \tilde z^o_{(i_W)}$. This can be done as follows. Let be $R>0$ such that $W\subset {\sf B}_R$. Simulate the restricted Poisson point process $\hat X^* \cap ([ {\sf B}_{nR}] \times (0, \frac{1}{n}))$, starting with $n=2$. If the generated zero-cell contains $W$, then choose it as $\tilde z^o_{(i_W)}$; otherwise update $n:=n+1$ until an appropriate zero-cell is simulated.

Then launch the (L-$G$)  (D-$\Lambda$) cell division process inside
$\tilde z^o_{(i_W)}$. Its restriction to $W$ has exactly the distribution of $(T_{t+\tau_{i_W}} \land W, t>0)$. Here,  the jump time $\tau_{i_W}$ of the zero-cell process is unknown because we do not know the value of the series \eqref{eq: def clock}.

\section{The Mondrian directional distribution}\label{sec:Mondrian}

Now we consider a particular class of directional distributions $\varphi$ and functionals $G$ where the sufficient condition of Theorem \ref{thm: sufficient existence} is satisfied.

Denote by $e_1,\ldots , e_d$ the orthonormal base of $\RR^d$ with 
$e_k := (0, \ldots ,0, 1,0\ldots 0)$ where the 1 is at the $k$-th position, $k\in \{1,\ldots , d\}$. Let $\delta_{e_k}$ denote the Dirac probability measure defined on $S^{d-1}$ and with mass 1 on $e_k$.

In this section we consider spherical directional distributions of the form   
\begin{equation}\label{cuboidal directional distr}
\varphi = \frac{1}{2} \sum_{k=1}^d p_k (\delta_{e_k}+\delta_{-e_k}) ,\quad p_k >0, \ \sum_{k=1}^d p_k =1.
\end{equation}
Referring to recent papers like \cite{OReilly_Tran2022}, we will call them Mondrian directional distributions.

The random tessellations generated by Mondrian directional distributions consist of cuboids. The following Theorems \ref{thm:existence volume weighted} and \ref{thm: intrinsic volumes} provide  first concrete models of  cell-division processes in $\RR^d$ that differ from the STIT model.

\begin{figure}\label{fig:stitFigures}
\centerline{
\includegraphics[width=0.450\textwidth]{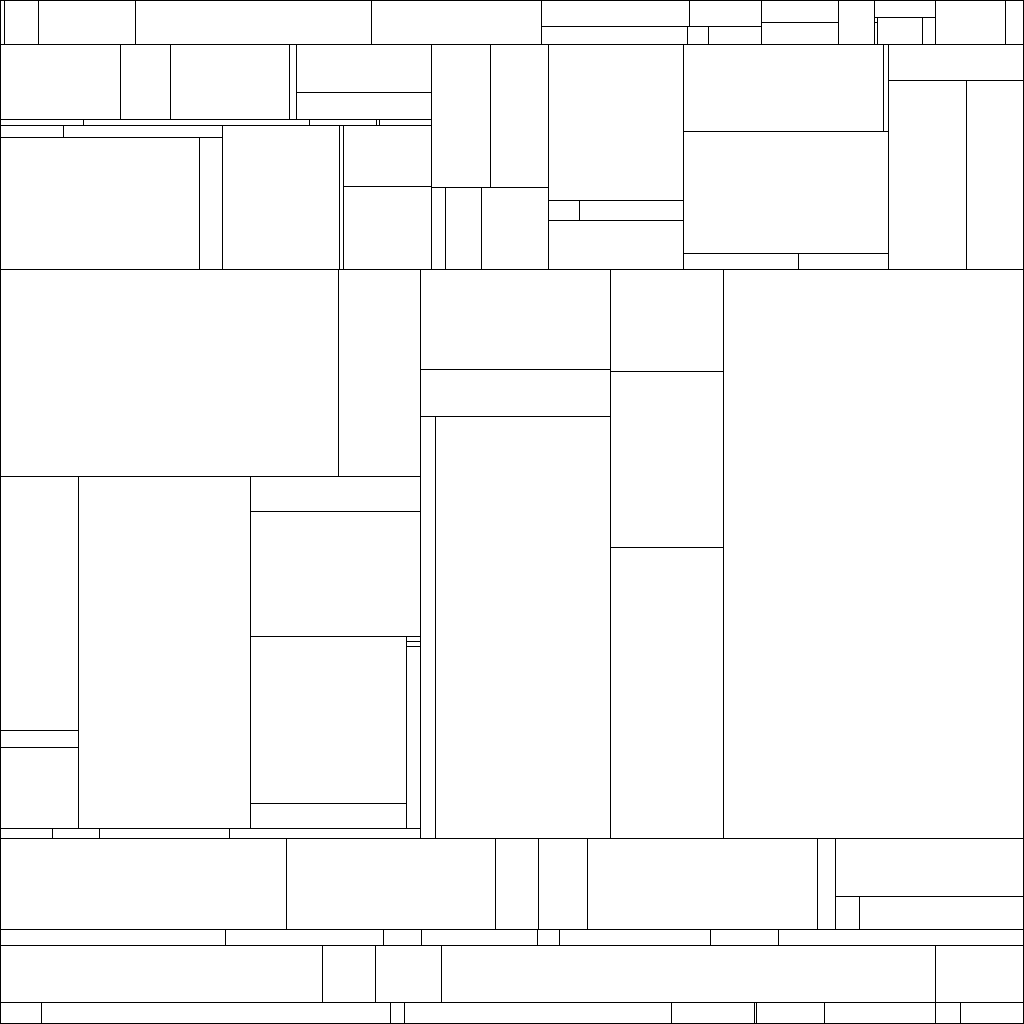}
\hspace{.5cm}
\includegraphics[width=0.450\textwidth]{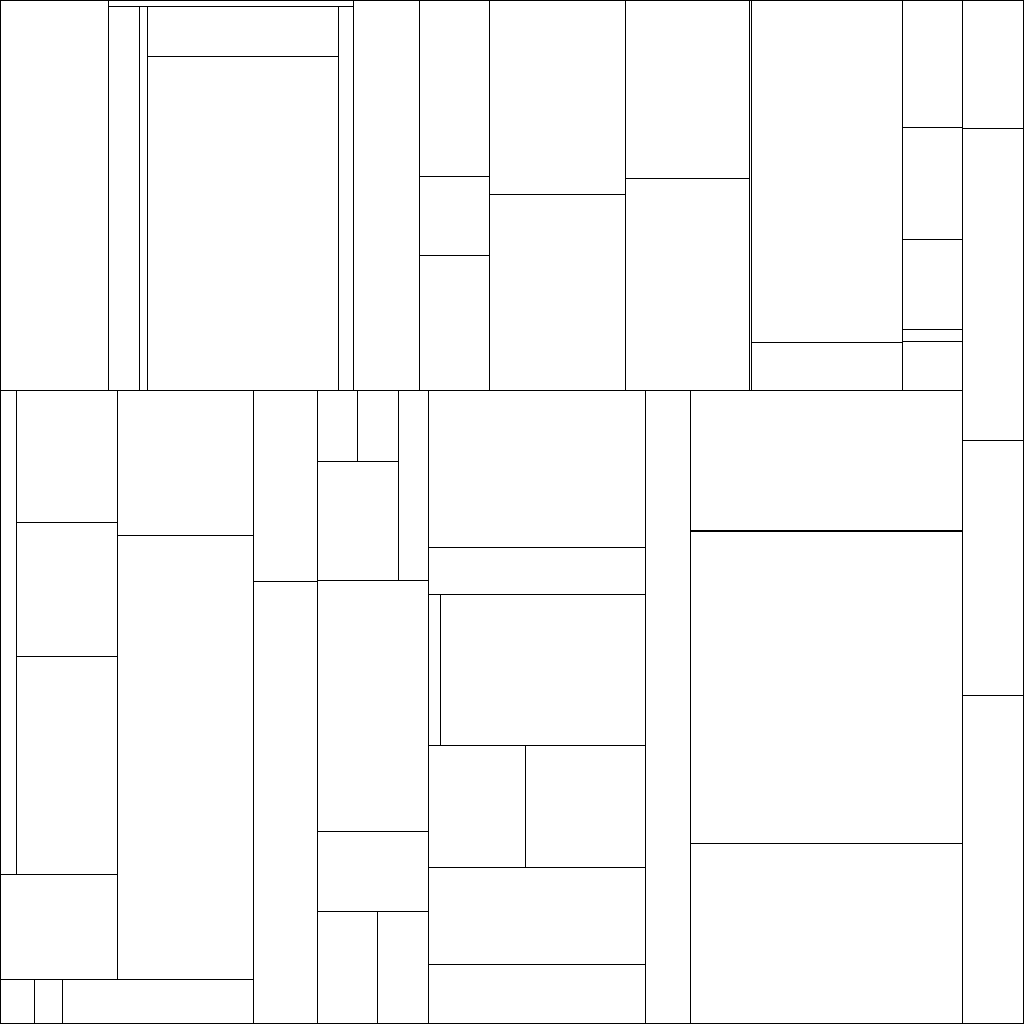}
}
\caption{Simulations of cell division tessellations in a quadratic window with (L-$\Lambda$) (left panel) and (L-$V_2$) (right panel), both with (D-$\Lambda$). The directional distribution is $\varphi =  \frac{1}{2} \delta_{e_1} + \frac{1}{2} \delta_{e_2}$. (Generated with the software \cite{crackpatternsimul}.)}
\end{figure}

\subsection{Existence results}

For a $d$-dimensional cuboid $z$ with all $(d-1)$-dimensional faces parallel or orthogonal to the axes in $\RR^d$, denote by by $l^{(1)},\ldots ,l^{(d)} $ the side lengths of the cuboid, where $l^{(k)}$ is the length of the sides that are parallel to $e_k$, $k\in \{ 1,\ldots , d\}$.  Denote $S(z):= \sum_{k=1}^d l^{(k)}$ if $z$ is a cuboid and $S(z)=1$ elsewhere.

\begin{theorem}\label{thm:existence volume weighted}
Let $\Lambda$ be a translation invariant measure on the space of hyperplanes in $\RR^d$, satisfying \eqref{eq:hyperplane measure}  where $\varphi$ is a Mondrian directional distribution \eqref{cuboidal directional distr}.  Then there exists a Markov tessellation-valued process $(T_t, t>0)$, which is an (L-$S$) (D-$\Lambda$) cell division process according to Definition \ref{def:cell div proc}.
\end{theorem}

{\bf Proof}
Consider the auxiliary Poisson process $\hat X^*$ of hyperplanes that are marked with birth times, that is the Poisson point process on $\hH \times (0,\infty )$ with the intensity measure $\Lambda \otimes \lambda_+$, see Section \ref{section:Notation}, and $\Lambda$ defined in \eqref{eq:hyperplane measure} with $\varphi$ as in \eqref{cuboidal directional distr}. We show that the sufficient condition of Theorem \ref{thm: sufficient existence} is satisfied for the  functional $G=S$, which is monotone and translation invariant on the class of cuboids. 

The Poisson point process $\hat X^*$ can be described as a superposition of $d$ i.i.d. Poisson point processes $\hat X^{(k)*}$, that is $\hat X^*\stackrel{D}{=}\bigcup_{k=1}^d \hat X^{(k)*}$, and $\hat X^{(k)*}$ is a Poisson point process on $\hH \times (0,\infty )$ with the intensity measure $p_k\, \Lambda^{(k)} \otimes \lambda_+$, and $\Lambda^{(k)}$ as defined in \eqref{eq:hyperplane measure} with $\varphi^{(k)} =\frac{1}{2}(\delta_{e_k}+ \delta_{-e_k})  $, $k\in \{ 1,\ldots , d\} $ (for which condition \eqref{assumption:hyperplane measure} is obviously not satisfied).

For $k\in \{ 1,\ldots , d\} $ we define the Poisson point processes
$\Phi^{k+}$ and $\Phi^{k-}$   on $(0,\infty )\times (0,\infty )$ as the birth time marked processes of intersection points of the hyperplanes of $\hat X^{(k)*}$ with the positive axis and the negative axis in direction  $e_k$ and $-e_k$, respectively. Formally,

$\Phi ^{k+}:= \{ (x,t) \in (0,\infty )\times (0,\infty ): (h(e_k, x),t) \in \hat X^{(k)*}\}$ 

and 

$\Phi ^{k-}:= \{ (x,t) \in (0,\infty )\times (0,\infty ): (h(-e_k, x),t) \in \hat X^{(k)*}\}$.

Both point processes have the intensity measure  $p_k\, \lambda_+ \otimes \lambda_+$.
Now we define the Markov chains

$M^{k+}=((x^{k+}_n,t^{k+}_n),\, n\in \NN_0)$

with $(x^{k+}_n,t^{k+}_n)\in \Phi^{k+}$ and

$x^{k+}_0:= \min \{ x:\, (x,t)\in \Phi^{k+} , t\leq 1 \}$,

$x^{k+}_{n+1}:= \min \{ x:\, (x,t)\in \Phi^{k+} , t\leq t^{k+}_n \}$, $n\in \NN_0$.

The times $t^{k+}_n$ are  a.s. uniquely defined by the condition $(x^{k+}_n,t^{k+}_n)\in \Phi^{k+}$.  Correspondingly,  $M^{k-}$ is defined, replacing $k+$ by $k-$.

The properties of the Poisson process $\hat X^*$   yield  that 
for fixed $k\in \{ 1,\ldots , d\} $ the Markov chains $M^{k+}$ and $M^{k-}$ are i.i.d.

Let us consider some properties of $M^{k+}$.
The random variable $x^{k+}_0$ is exponentially distributed with parameter $p_k$, and $t^{k+}_0$ is uniformly distributed in the interval $(0,1)$. Then, $x^{k+}_{n+1}-x^{k+}_n$ is exponentially distributed with the parameter $p_k\, t^{k+}_n$ and independent of $x^{k+}_0,\ldots , x^{k+}_n$ and, given $t^{k+}_n$,
conditionally independent of  $t^{k+}_0,\ldots , t^{k+}_{n-1}$. Furthermore, $t^{k+}_{n+1}$ is uniformly distributed in the interval $(0,t^{k+}_n)$ and, given $t^{k+}_n$,
conditionally independent of $x^{k+}_0,\ldots , x^{k+}_n$ and of $t^{k+}_0,\ldots , t^{k+}_{n-1}$.

Now consider the Markov chain $((S_{i},t_{i}),\, i \in \ZZ_-)$ where
$(t_i ,\, i \in \ZZ )$ is the ordered sequence of jump times of $(\tilde z^o_t , t>0)$ with $t_0<1<t_1$, as defined in Section \ref{section:Notation}, and $S_i:= S\left(\tilde z^o_{(i)}\right)$.
Thus we have 
$t_0 = \max \{ t^{k+}_0,  t^{k-}_0 , k\in \{ 1,\ldots , d\}  \}$, and 
\begin{equation*}\label{eq:sequence S initial}
S_0= \sum_{k=1}^d  x^{k+}_0 + x^{k-}_0 . 
\end{equation*}
By recursion, 
$
t_{i-1}=  \max \{ t^{k+}_n <t_i,  t^{k-}_n<t_i , k\in \{ 1,\ldots , d\} , \ n\in \ZZ_- \},
$
and 
\begin{equation}\label{eq:sequence S}
S_{i-1} = \begin{cases}S_i + x^{k+}_{m+1}-x^{k+}_m  & \textup{ if } t_{i} =t^{k+}_m \ \textup{ for some } m\in \ZZ_- \ ,\\
	                                     S_i + x^{k-}_{m+1}-x^{k-}_m  &  \textup{ if } t_i =t^{k-}_m \ \textup{ for some } m\in \ZZ_- \ .
													\end {cases}						
\end{equation}

Let $(R_i, \, i\in \ZZ )$ be a sequence of i.i.d. random variables with the probability density $f(r)= 2d r^{2d-1} {\bf 1}\{ 0<r<1 \}$, that is, $R_i$ has the distribution of the maximum of $2d$\,  i.i.d. random variables that are uniformly distributed in the interval $(0,1)$. 

Thus the law of  the time component $t_i$ in the Markov chain $((S_{i},t_{i}),\, i \in \ZZ_-)$ can be described by
$$
t_0 \stackrel{D}{=} R_0\, ,
$$
and by recursion
$$
t_{i-1} \stackrel{D}{=} t_i \,  R_{i-1}\, ,\quad i \in \ZZ_-\, .
$$
Thus
$$
(t_i,\, \, i \in \ZZ_- ) \stackrel{D}{=} \textstyle{\left(\prod_{j=i}^{0} R_j,\, \, i \in \ZZ_- \right)} \, .
$$

The  law of $S_{i-1} -S_i$ depends on both $t_{i}$ and the index $k^+$ or $k^-$.  It is the exponential distribution with the parameter $p_k \, t_{i}$, not depending on the sign of the index $k$.

Let $(\kappa_i ,\, i\in \ZZ )$ be a sequence of i.i.d. random variables, exponentially distributed with parameter 1. This sequence is assumed to be independent of all the other random variables we considered so far.
Then the conditional distribution is
\begin{align}\label{eq:length in direction k} 
& S_{i-1} -S_i \stackrel{D}{=}    \left( p_k \, \prod_{j=i}^0 R_j\right)^{-1}\kappa_i \leq p_{\min}^{-1} \left(  \prod_{j=i}^0 R_j\right)^{-1}\kappa_i, \\ 
&\textup{under the condition that  } t_i =t^{k+}_m \textup{ or } t_i =t^{k-}_m \textup{ for some } m \in \NN_0\, , \nonumber
\end{align}
with $p_{\min} := \min \{ p_1,\ldots , p_k\}$.

Thus a sufficient condition for 
$$
\sum_{i\leq 0} S_i^{-1} = S_0^{-1} + \sum_{i\leq 0}\left(S_0 + \sum_{\ell=i}^0 (S_{\ell-1} -S_\ell) \right)^{-1}
<\infty \quad  a.s.
$$
is
\begin{equation}\label{eq:sufficient cond cuboidal}
\sum_{i\leq 0} \left(\sum_{\ell=i}^0\left(  \prod_{j=\ell}^0 R_j\right)^{-1} \kappa_\ell \right)^{-1}  <\infty \quad  a.s.
\end{equation}

Note that in \eqref{eq:sufficient cond cuboidal} the item $S_0^{-1} $ is neglected in order to simplify the technicalities.
To verify \eqref{eq:sufficient cond cuboidal} we provide an upper bound for the expectation of the random variable on the left-hand side.

\begin{align*}
& \EE \sum_{i\leq 0} \left(\sum_{\ell=i}^0\left( \prod_{j=\ell}^0 R_j\right)^{-1} \kappa_\ell\right)^{-1} \! \! \!  \! \! \! \!  \! \! \! \!  \!
&\!\leq &\ \EE \sum_{i\leq 0} \left(\sum_{\ell=i}^{i+1}\left(  \prod_{j=\ell}^0 R_j\right)^{-1} \kappa_\ell\right)^{-1} \\
& & &\\
= & \, \EE \sum_{i\leq 0}
    \frac{\prod_{j=i+1}^{0} R_j}{ R_i^{-1} \, \kappa_i +\kappa_{i+1} } 
&\leq & \, \EE \sum_{i\leq 0} 
    \frac{\prod_{j=i+1}^{0} R_j}{ \kappa_i + \kappa_{i+1}  } \\
& &&\\
= & \, \EE\frac{1}{ \kappa_0 + \kappa_{1}} \ \sum_{i\leq 0}  \left(\prod_{j=i+1}^{0} \EE R_j\right) \, 
&= & \,  \sum_{i\leq 0}  \left(\prod_{j=i+1}^{0} \frac{2d}{2d+1}\right) 
< \infty \, .
\end{align*}
In the last  equation we used that the sum of two i.i.d. exponentially distributed random variables with parameter 1, has a gamma-distribution with parameter $(2,1)$, also referred to as an Erlang distribution, and hence $\EE(\kappa_0 + \kappa_{1})^{-1}=1$. Thus it is shown that \eqref{eq:sufficient cond cuboidal} is satisfied.
\hfill $\Box$

In order to generalize the result for Mondrian directional distributions, we consider the individual side lengths of the zero-cells. 
Denote by $l_i^{(1)},\ldots ,l_i^{(d)} $ the side lengths of the cuboid $\tilde z^o_{(i)}$, where $l_i^{(k)}$ is the length of the sides that are parallel to $e_k$, $k\in \{ 1,\ldots , d\}$, $i\in \ZZ$. The following technical lemma  prepares the proof of Theorem \ref{thm: intrinsic volumes}.
\begin{lemma}\label{lemma: all sides longer than 1}
With the assumptions of Therorem \ref{thm:existence volume weighted},
$$
\PP \left(\exists i_0\in \ZZ \, \forall i\leq i_0 \, \forall k\in \{ 1,\ldots , d\}:\, l_i^{(k)}>d \right) =1.
$$
\end{lemma}
{\bf Proof}
Fix some  $k\in \{1,\ldots , d\}$ and, using the notation in the proof of Theorem \ref{thm:existence volume weighted}, consider \eqref{eq:length in direction k} which  expresses the growth of the length $l_i^{(k)}$ if $ t_i =t^{k+}_m \textup{ or } t_i =t^{k-}_m \textup{ for some } m \in \NN_0$.
The index $k$ is chosen with probability $p_k$. Obviously, \\
 $\left( p_k \, \prod_{j=i}^0 R_j\right)^{-1} > 1$.
Hence, the  events 
$$
B_i:= \{ \kappa_i >d,\,  t_i =t^{k+}_m \textup{ or } t_i =t^{k-}_m \textup{ for some } m \in \NN_0\}, \quad i\in \ZZ,
$$
which are independent, have a strictly positive probability that do not depend on $i$. Thus the Borel-Cantelli Lemma yields that
$\PP \left(  \bigcap_{n\leq 0} \bigcup_{i\leq n} B_i \right) =1 ,$
that is there occur infinitely many of the $B_i$ a.s. 
This, together with the monotonicity of the sequence $(l_i^{(k)}, \, i\in \ZZ)$ yields that for all $k\in \{1,\ldots , d\}$ we have
$$
\PP \left(\exists i_0(k)\in \ZZ \, \forall i\leq i_0(k) :\, l_i^{(k)}>d\right) =1,
$$
and hence 
$$
\PP \left(\forall k\in \{1,\ldots , d\}\, \exists i_0(k)\in \ZZ \, \forall i\leq i_0(k) :\, l_i^{(k)}>d\right) =1.
$$
Because there are only finitely many $k$, the assertion of the lemma follows.
\hfill $\Box $

Lemma \ref{lemma: all sides longer than 1} can be used to check the sufficient condition \eqref{eq: suff cond} when the functional $G$ is chosen as a function of the side lengths of the zero cell. Here we consider the important cases of the intrinsic volumes $V_n$, $n\in \{ 1,\ldots ,d \}$. The functional  $V_d$ is the volume, $V_1$ the mean width (or breadth) up to a constant factor. For $n=3$ the intrinsic volume $V_2$ is, up to a constant factor, the surface area. The intrinsic volume $V_0$ is the Euler characteristic which has the value 1 for all nonempty convex sets, and hence it is clear that for $G=V_0$ the sufficient condition \eqref{eq: suff cond} is not satisfied. 

A definition of the intrinsic volumes for convex bodies can be found in \cite{sw08}. Here we make use of a formula for cuboids which is a particular case of the formula for polytopes, see (14.35) of \cite{sw08}. For $n\geq 1$
\begin{equation}\label{eq: intrinsic volumes for cuboids}
V_n (\tilde z^o_{(i)})=  \sum_{1\leq k_1<\ldots <k_n \leq d} \ \ 
 \prod_{j=1}^n l_i^{(k_j)} \, .
\end{equation}

\begin{theorem}\label{thm: intrinsic volumes}
Let $\Lambda$ be a translation invariant measure on the space of hyperplanes in $\RR^d$, satisfying \eqref{eq:hyperplane measure}  where $\varphi$ is a Mondrian directional distribution \eqref{cuboidal directional distr}. Then, for all $n\in \{ 1,\ldots ,d \}$, there exists a Markov tessellation-valued process $(T_t, t>0)$, which is an (L-$V_n$) (D-$\Lambda$) cell division process according to Definition \ref{def:cell div proc}.
\end{theorem}
{\bf Proof}
It is sufficient to show that \eqref{eq: suff cond} is satisfied for $G=V_n$ for all  $n\in \{ 1,\ldots ,d \}$.
For $n=1$ it follows immediately from Theorem \ref{thm:existence volume weighted}.
For $n \in\{2,\ldots ,k\}$ and $i\in \ZZ$,
$$
\left(\sum_{k_1<\ldots <k_n} \, \prod_{j=1}^n l_i^{(k_j)}\right)  - \sum_{k=1}^d  l_i^{(k)}=
\frac{1}{{d-1 \choose n -1}} \sum_{k=1}^d  l_i^{(k)} \left(  \sum_{{k_1<\ldots <k_{n-1}}\atop {k_j\not= k}} \left(\frac{{d\choose n}}{d}\prod_{j=1}^{n -1} l_i^{(k_j)} -1 \right) \right) .
$$
Note that ${d \choose d}/d =1/d$ and ${d \choose n}/d \geq 1$ for $n\in \{ 1,\ldots , d-1\}$.
By Lemma \ref{lemma: all sides longer than 1} there exists a.s. an $i_0\in \ZZ$ such that for all $i\leq i_0$ and all $k\in \{1,\ldots , d\}$ the side lengths $ l_i^{(k)} > d$, which implies that the difference on the left-hand side is positive.
Recall that $S_i=\sum_{k=1}^d  l_i^{(k)}$ for all $i\in \ZZ$, and that in the proof of  of Theorem \ref{thm:existence volume weighted} it was shown that 
$
\sum_{i\leq 0} S_i^{-1} 
<\infty \quad  a.s.
$
Hence, for $n \in\{2,\ldots ,k\}$ also
$$
\sum_{i\leq 0} \left(\sum_{k_1<\ldots <k_n} \, \prod_{j=1}^n l_i^{(k_j)}\right)^{-1}  <\infty \quad a.s.
$$
and  by \eqref{eq: intrinsic volumes for cuboids},
$$
\sum_{i\leq 0} V_n\left(\tilde z^o_{(i)}\right)^{-1} \, 
                < \infty  \quad a.s.
$$
\hfill $\Box $

\begin{remark}\label{rem:exponent alpha}
Theorem \ref{thm: intrinsic volumes} can easily be generalized to 
(L-$V_n^\alpha$) (D-$\Lambda$) cell division process, for exponents $\alpha \geq 1$. The larger $\alpha$, the longer will be the life time of cells $z$ with $V_n(z) <1$ and the shorter will be the life time if $V_n(z)>1$. 
\end{remark}

Note that for the Mondrian distribution defined in \eqref{cuboidal directional distr} with $p_1=\ldots =p_d= 1/d$, the (L-$V_1$) (D-$\Lambda$) cell division process is -- up to a scaling factor -- the same as the STIT tessellation process driven by $\Lambda$, which is a (L-$\Lambda$) (D-$\Lambda$) cell division process.

\subsection{ Distribution of cell volumes} 

Now we consider the particular cell division process in $\RR^d$ 
with volume-weighted life time distribution (L-$V_d$) and (D-$\Lambda $),   with a Mondrian directional distribution \eqref{cuboidal directional distr}. We make use of an idea used by Cowan in \cite{cowan2010} developed for so-called 'geometry-independent apportionment of volume' (GIA)  which means a division rule such that the ratio of the volume of a daughter cell and the volume of the mother cell is uniformly distributed on the interval $(0,1)$.

In contrast to this, in Subsection \ref{subsec:fragmentation} we will consider an (L-$V_d$) (D-$\Lambda $) process in a {\em bounded} window $W$ and describe its relation to fragmentation.

For a random stationary tessellation of $\RR^d$ the concept of the 'typical cell' is formally defined using the Palm measure,  see p. 450 of \cite{sw08} for example.
Intuitively, it can be imagined as a 'cell, cosen at random, where all cells have the same chance to be chosen'.

\begin{lemma} \label{lem: Poiss point on line}
Let $\Lambda$ be a translation invariant measure on the space of hyperplanes in $\RR^d$, satisfying \eqref{eq:hyperplane measure}  where $\varphi$ is a Mondrian directional distribution \eqref{cuboidal directional distr}, and let $(T_t,\, t>0)$ be the (L-$V_d$) (D-$\Lambda $) cell division process  as described in the proof of Theorem \ref{thm: sufficient existence}.
Then for all $t>0$ the law of the volume $V_d(z)$ of the typical cell of $T_t$ is the exponential distribution with parameter $t$, and the law of the volume $V_d(z^o_t)$ of the zero-cell of $T_t$ is the gamma-distribution with parameter $(2,t)$, which is also referred to as an Erlang distribution.
\end{lemma}
{\bf Proof}
The volumes of cells of a tessellation $T_t$, $t>0$,  are represented by the lengths of intervals on $\RR$.
The idea of the proof is the following: The tessellation $T_t$ is mapped to a point process $(y^{(j)},\, j \in \ZZ)$ on $\RR$ such that the lengths of the intervals between adjacent points correspond to the volumes of the cells of $T_t$.

In order to formalize this idea, it is an involved  technical issue  to bring these intervals into an appropriate linear order.
As before, we start with the process of zero cells as defined in  \eqref{eq: def zero} and we use a method analogous to that one in the proof of Theorem \ref{thm:existence volume weighted}. But here we consider the differences  of volumes when cells are split.

If for some $\eps >0$ the tessellation $T_\eps$ is mapped to a point process $(y_\eps ^{(j)},\, j \in \ZZ)$ on $\RR$ such that the lengths of the intervals between adjacent points correspond to the volumes of the cells of $T_\eps$, then the division of the cell volumes in the continuation of the cell division process corresponds to a {\em Poisson} point process that is superposed to $(y_\eps ^{(j)},\, j \in \ZZ)$. Because this can be done for any $\eps >0$, one can conclude that the cell volumes of $T_t$ are represented by the lengths of the interval between two adjacent points of a Poisson point process.
Summarizing the idea, we will show that for a Poisson point process on $\RR^d \times (0,\infty )$ with intensity maesure $\lambda \otimes \lambda_+$

Formalizing this idea, it is an issue  to bring these intervals into an appropriate linear order.
Let us begin with the process of zero cells as defined in  \eqref{eq: def zero}. We use a method analogous to that one in the proof of Theorem \ref{thm:existence volume weighted} above, but here we consider the differences  of volumes when cells are split.

Fix an $\eps >0$ and $i\in \ZZ$ such that $\tau_i\leq \eps < \tau_{i+1}$, where $\tau_i$ is defined in \eqref{eq: def clock} with $G=V_d$.
Let  $\hat y_\eps^{(-1)} < 0 < \hat y_\eps^{(0)}$ be real valued such that 
$ \hat y_\eps^{(0)} - \hat y_\eps^{(-1)}= V_d (z_{\tau_i})$, and $0$ is uniformly distributed in the interval $(\hat y_\eps^{(-1)}, \hat y_\eps^{(0)})$ 

Initialize a counter $r:=0$ and a counter $\ell :=1$. 
If in \eqref{eq:sequence S} $t_i =t^{k+}_m$ for some $k$ and $m$, then  put
$\hat y_\eps^{(1)}:=\hat y_\eps^{(0)} + V_d(z_{\tau_{i-1}}\setminus z_{\tau_i})$, $t^{(0)} := \tau_i$, and update $r:=1$.
If  $t_i =t^{k-}_m$ for some $k$ and $m$, then  define
$\hat y_\eps^{(-2)}:=\hat y_\eps^{(-1)} -V_d(z_{\tau_{i-1}}\setminus z_{\tau_i})$, $t^{(-1)} := \tau_i$, and update $\ell:=2$.

Having defined $(\hat y_\eps^{(-\ell)},t^{(-\ell)}),\ldots , (\hat y_\eps^{(r)},t^{(r)})$ for some $\ell \in \NN_0$ and $r\in \NN_0$, 
the next item of the sequence is:
If in \eqref{eq:sequence S} $t_{i-r-\ell -1} =t^{k+}_m$ for some $k$ and $m$, then
$\hat y_\eps^{(r+1)}:=\hat y_\eps^{(r)} + V_d(z_{\tau_{i-r-\ell -2}}\setminus z_{\tau_{i-r-\ell -1}})$, 
$t^{(r+1)}:=\tau_{i-r-\ell -1}$,
and update $r:=r+1$.
If  $t_{i-r-\ell -1} =t^{k-}_m$ for some $k$ and $m$, then
$\hat y_\eps^{(-\ell -1))}:=\hat y_\eps^{(-\ell )}-V_d(z_{\tau_{i-r-\ell -2}}\setminus z_{\tau_{i-r-\ell -1}})$, $t^{(-\ell -1)}:=\tau_{i-r-\ell -1}$, and update $\ell :=\ell +1$.

This yields the sequence $((\hat y_\eps^{(j)},t^{(j)}) ,\, j\in \ZZ)$ pertaining to the process of zero cells.
It has to be complemented by the points which correspond to the divisions of the cells $\textup{cl}(z_{\tau_{i-r-\ell -2}}\setminus z_{\tau_{i-r-\ell -1}})$ in the time interval $(\tau_{i-r-\ell -1}, t)$.

 This can be described as follows.
The volume of the cell  $\textup{cl}(z_{\tau_{i-r-\ell -2}}\setminus z_{\tau_{i-r-\ell -1}})$ equals the length of the interval $(\hat y_\eps^{(r)}, \hat y_\eps^{(r+1)})$ or $(\hat y_\eps^{(-\ell -1)},\hat y_\eps^{(-\ell )})$, respectively. For the sake of simplicity, in the following we consider only the part of $((\hat y_\eps^{(j)},t^{(j)}) ,\, j\in \ZZ)$ which belongs to $(0,\infty )\times (0,\infty )$. The other part on 
$(-\infty ,0)\times (0,\infty )$  can be dealt with quite analogously.
When the cell  $\textup{cl}(z_{\tau_{i-r-\ell -2}}\setminus z_{\tau_{i-r-\ell -1}})$ is divided by a hyperplane, the respective interval is divided by a point such that the two new interval lengths equal the two volumes of the  daughter cells. We define the new interval which is closer to $0\in \RR$ to pertain to the daughter cell which is contained in the half space of the dividing hyperplane that contains the origin $o\in  \RR^d$. The point dividing the interval is marked with the birth time of the hyperplane dividing the cell.
Thus the interval $(\hat y_\eps^{(r)}, \hat y_\eps^{(r+1)})$ is divided after a life time, which is exponentially distributed with the parameter  $\hat y_\eps^{(r+1)}- \hat y_\eps^{(r)}$, and the dividing point is -- due to the translation invariance of $\Lambda$ -- uniformly distributed in this interval. The  new daughter intervals are subsequently divided, following analogous rules. Thus the birth time marked point process appearing in the interval  $(\hat y_\eps^{(r)}, \hat y_\eps^{(r+1)})$ has the same distribution as the Poisson point process on $(\hat y_\eps^{(r)}, \hat y_\eps^{(r+1)})\times (\tau_{i-r-\ell -1}, t)$ with the two-dimensional  Lebesgue measure (restricted to this set) as its intensity measure.

Taking all the marked points together, we obtain a point process\\ $Y_\eps := ((y_\eps^{(j)}, t^{(j)}),\, j\in \ZZ )$. 

Given the point process $Y_\eps = ((y_\eps^{(j)}, t^{(j)}),\, j\in \ZZ )$ for some time $\eps >0$, and regarding the independence assumptions in Definition \ref{def:cell div proc}  the generation of  the point process $Y_t =((y_t^{(j)}, t^{(j)}),\, j\in \ZZ )$ for $t>\eps $ can be described 
using an auxiliary Poisson point process $\Phi^*$   on $\RR \times (0,\infty )$ with intensity measure $\lambda \otimes \lambda_+$. 

For all $0<\eps <t$ let be $\Phi^*_{(\eps ,t)} := \{ (x,t')\in \Phi^* : \,  \eps  <t' <t \}$. We obtain that
\begin{equation}\label{eq: point process volumes}
Y_t \stackrel{D}{=} Y_\eps \cup \Phi^*_{(\eps ,t)} ,
\end{equation}
which is the superposition of two independent point processes.

The sequence $((\hat y_\eps^{(j)},t^{(j)}) ,\, j\in \ZZ)$ pertaining to the process of zero cells is monotone in $\eps$ in the following sense. For $\eps_2 < \eps_1$ the sequence $((\hat y_{\eps_2}^{(j)},t^{(j)}) ,\, j\in \ZZ)$ emerges from $((\hat y_{\eps_1}^{(i)},t^{(i)}) ,\, i\in \ZZ)$ by deleting the points with
 $$
t^{(i+1)} >\eps_2 \mbox{ for } i\geq 0,\qquad
 \mbox{or}\qquad
 t^{(i-1)} >\eps_2 \mbox{ for } i< 0, 
$$
and after this, a respective rearrangement of the indexes $i\in \ZZ$.\\
Accordingly, $Y_{\eps_2}=((y_{\eps_2}^{(j)}, t^{(j)}),\, j\in \ZZ )$ can -- up to a rearrangement of the indexes -- be considered as a subsequence of $Y_{\eps_1}=((y_{\eps_1}^{(j)}, t^{(j)}),\, j\in \ZZ )$. 

Denoting by $Y_{(\eps_2 ,\eps_1)}$ the difference between $Y_{\eps_1}$ and $Y_{\eps_2}$  we can  write 
\begin{equation*}
Y_{\eps_2} \cup \Phi^*_{(\eps_2 ,t)}  \stackrel{D}{=} Y_{\eps_2} \cup \Phi^*_{(\eps_2 ,\eps_1)}   \cup \Phi^*_{(\eps_1 ,t)} \stackrel{D}{=} Y_{\eps_2}
\cup Y_{(\eps_2 ,\eps_1)}
 \cup \Phi^*_{(\eps_1 ,t)} 
,
\end{equation*}
where the superposed point processes are independent.
This yields that $ Y_{(\eps_2 ,\eps_1)} \stackrel{D}{=} \Phi^*_{(\eps_2 ,\eps_1)}$.  
Intuitively, this means that the sequence $((\hat y_{\eps_1}^{(i)},t^{(i)}) ,\, i\in \ZZ)$ pertaining to the process of zero cells until time $\eps_1$ is consistent with the Poisson point process properties, which are characteristic for the consecutive division of cell volumes.

Because this holds for all $\eps_2 > \eps_1 >0$, it is shown that $Y_t=((y_t^{(j)}, t^{(j)}),\, j\in \ZZ )$ is a Poisson point process on $\RR \times (0,t)$ with intensity measure $\lambda \otimes \lambda_{(0,t)}$, where $\lambda_{(0,t)}$ denotes the restriction of the Lebesgue measure to the interval $(0,t)$. Hence the point process $(y_t^{(j)},\, j\in \ZZ )$ of  points projected onto $\RR$, is a homogeneous Poisson point process with intensity measure $t\, \lambda$, and therefore the length of the typical interval between two points is exponentially distributed with parameter $t$, while the length of the interval containing the origin $o$ has a gamma-distribution with parameter $(2,t)$.
\hfill $\Box$

Cell division processes are potential models for real crack or fracture pattern appearing in materials science, geology, soft matter or nanotechnology. The examples considered in \cite{nagel_mecke_ohser_weiss2008, leon_ohser_nagel_arscott2019} illustrate, that a useful and simple initial goodness-of-fit criterion for such models is the coefficient of variation CV of the volume of the typical cell. And it became obvious that this CV of STIT is too large compared to the CV for data of real crack pattern.

\begin{corollary}
With the assumptions of Lemma \ref{lem: Poiss point on line} for all $t>0$
$$\textup{CV}(V_d (z)):= \frac{\sqrt{\textup{Var} (V_d (z))}}{\EE V_d (z)}=1,$$
where $z$ is a random cuboid with the distribution of the typical cell of $T_t$ and $\textup{Var}$ denotes the variance of a random variable.
\end{corollary}

For the cell division process with (L-$\Lambda$) (D-$\Lambda$), which is the STIT tessellation process, at any time $t>0$ the typical cell has the same distribution as typical cell of Poisson hyperplane tessellation for the same $\Lambda$ at the same time, see \cite{nagel_weiss2005}.
Hence, for the translation invariant measure $\Lambda$   satisfying \eqref{eq:hyperplane measure}  where $\varphi$ is a Mondrian directional distribution \eqref{cuboidal directional distr}, the CV of the volume of the typical cell equals  $\textup{CV}\left(\prod_{j=1}^d Z_j\right)$, where $Z_1,\ldots ,Z_d$ are independent exponentially distributed random variables with the parameters $p_1 t,\ldots ,p_d t$, respectively.
Thus 
$$\textup{CV}\left(V_d \left(z^{\textup{STIT}}\right)\right)=2^d -1,$$
where $z^{\textup{STIT}}$ is a random cuboid with the distribution of the typical cell of the STIT tessellation.

Note that these CV do neither depend on the time $t$ nor on the probabilities $p_1 ,\ldots ,p_d $.
Obviously, for $d\geq 2$ the variability of the volume of the typical cell of the (L-$V_d$) (D-$\Lambda $) cell division process is considerably smaller than that one of STIT.

\subsection{A cell division process in a bounded window $W$ and its relation to fragmentation}\label{subsec:fragmentation}

Now consider a fixed window $W\in \hP$ which is a cuboid with its sides parallel or orthogonal to all axes of $\RR^d$. An (L-$V_d$) (D-$\Lambda $) cell division process $(T_{W,t},\, t\geq 0)$ {\em within} $W$ is launched with $T_{W,0}:= \{W\}$. Note that its distribution is different from $(T_t \wedge W ,\, t>0)$, which is the restriction to $W$ of the (L-$V_d$) (D-$\Lambda $) cell division process$(T_t ,\, t>0)$  in the whole $\RR^d$. 

The process $(T_{W,t},\, t\geq 0)$ can be related to a particular case of a conservative fragmentation process as described in \cite{bertoin2006}.

Assume that $G(W) =1$. For $t>0$ and $T_{W,t} =\{ z_1,\ldots ,z_n\}$      define 
\begin{equation*}
{\bf G}(T_{W,t}):=( G(z_{i_1}),\ldots , G(z_{i_n}), 0,0,\ldots  )
\end{equation*}
such that $\{i_1,\ldots , i_n\} =\{ 1,\ldots ,n\}$ and $G(z_{i_1})\geq \ldots \geq G(z_{i_n})$.

If we choose $G(z)=V_d(z)$ for $z\in \hP$, then
 the (L-$V_d$) (D-$\Lambda $) cell division process $(T_{W,t},\, t\geq 0)$ induces a pure jump Markov process $({\bf G}(T_{W,t}),\, t\geq 0)$ on ${\mathfrak S}:=  (0,\infty )^\NN$ with the following properties.

\begin{enumerate}
\item[(a)] $\sum_{z \in T_{W,t}} G(z) =1$ for all $t\geq 0$.

\item[(b)] The holding time of a state $( G(z_{i_1}),\ldots , G(z_{i_n}),0,0,\ldots  )$ is exponentially distributed with parameter 1.

\item[(c)] At the $n$-th jump time $\tau_{n}$, $n\in \NN$, the process jumps from a state\\
 $( G(z_{i_1}),\ldots , G(z_{i_n}),0,0,\ldots    )$ to 
$( G(z_{j_1}),\ldots , G(z_{j_{n+1}}),0,0,\ldots    )$ by the following rule:\\
The index $i_k$, $k\in \{ 1,\ldots , n\}$, is chosen at random with probability $G(z_{i_k})= G(z_{i_k})/\sum_{\ell \in \{ 1,\ldots , n\}}G(z_{i_\ell})$. Then the cell $z_{i_k}$ is divided into two new fragments, $z'_{i_k}, z''_{i_k}$ say, and 
\begin{equation}\label{eq:disloc frag}
\left(G(z'_{i_k}), G(z''_{i_k})\right)=\left(U\, G(z_{i_k}), (1-U)\, G(z_{i_k})\right),
\end{equation}
where $U$ is a random variable, uniformly distributed on $(0,1)$ and independent of all the other random variables considered here.
\item[(d)] Putting the $G$-values of the $n+1$ cells into an descending order yields the new state $( G(z_{j_1}),\ldots , G(z_{j_{n+1}}),0,0,\ldots    )$. 
\end{enumerate}

This shows that, according to Definition 1.1 of \cite{bertoin2006}, the process $({\bf G}(T_{W,t}),\, t\geq 0)$ is a self-similar fragmentation chain with index of self-similarity $\alpha =1$ and the dislocation measure $\nu$ which is defined by \eqref{eq:disloc frag}, that is, $\nu$ is the law of
$
(\xi , 1-\xi ,0,0,\ldots ) ,
$
where $\xi$ is a random variable which is uniformly distributed on the interval $\left[\frac{1}{2}  , 1\right]$.

Note, if $G$ is chosen as another intrinsic volume $V_n$, $n<d$, then $\sum_{z \in T_{W,t}} G(z) \geq 1$ with $'>'$ after the first division. Then this is no longer the {\em conservative} scenario in the sense of \cite{bertoin2006}.

\section{Discussion and open problems}

The 'constructive proof' of the existence result of Theorem \ref{thm: sufficient existence} of a cell division process does not imply a uniqueness result. 
We conjecture that the distribution of the process is uniquely determined by $G$ and $\Lambda$, but it is an open problem, whether for (L-$G$) and (D-$\Lambda$) also other distributions of cell division processes exist. Up to now, only the uniqueness of the (L-$\Lambda$) (D-$\Lambda$) cell division process is shown, which is the STIT tessellation process driven by $\Lambda$.

Somehow related to the uniqueness problem is the question whether the sufficient condition \eqref{eq: suff cond} for the existence of an (L-$G$) (D-$\Lambda$) cell division process is also a  necessary one.

In a cell division process $(T_t,\, t>0)$ the cells $z$ are divided subsequently as long as $G(z)>0$. Therefore, in general there will not be a limiting random tessellation for $t\to \infty$. Thus the problem arises whether there is a scaling of the process which stabilizes the one-dimensional distributions of the process or even yields a time-stationary process.
For a tessellation $T\in \TT$ the scaling by a factor  $c>0$ means $cT:= \{ cz;\, z\in T\}$, with $cz:=\{ cx:\, x\in z\}$.

For the STIT tessellation process $(Y_t,\, t>0)$ driven by an arbitrary $\Lambda$ satisfying \eqref{eq:hyperplane measure} and \eqref{assumption:hyperplane measure} it was shown in \cite{nagel_weiss2005} that
$
t\, Y_t \stackrel{D}{=} Y_1 \ \textup{ for all } t>0.
$
Furthermore, in \cite{martinez_nagel2012} it was proven that for $a>1$ the process $a^t\, Y_{a^t}$ is stationary in time.

For (L-$V_d$) (D-$\Lambda$) and the Mondrian directional distribution, by Lemma \ref{lem: Poiss point on line}, for all $t>0$ and for the zero cell $z^0_t$ of $T_t$ we have 
$
t\, V_d (z^0_t)\stackrel{D}{=}V_d (z^0_1)
$
and hence
$
V_d (t^{1/d} z^0_t)\stackrel{D}{=}V_d (z^0_1).
$
This does not yet show a scaling property for the whole random tessellation $T_t$ but one can see that, if there is a scaling factor, it can only be $t^{1/d}$.

Regarding tail triviality and other ergodic properties of the tessellations generated by cell division, the method of the so-called encapsulation used in our earlier papers \cite{martinez_nagel2014, martinez_nagel20016}, cannot be applied here directly for the models which are not spatially consistent. But the method used in the proof of Lemma \ref{lemma:trivial tail zero cell} for zero-cells can be developed further. 

Real crack structures, some of them mentioned in the introduction, show the tendency that a crack of a cell divides it close to the center. This is not yet incorporated in the models with a (D-$\Lambda$) division rule where  $\Lambda$ is a translation invariant measure. 
Maybe that this can partially be compensated in some cases using the life time rule (L-$G^\alpha$) instead of (L-$G$), with a large value of the exponent $\alpha$, see also Remark \ref{rem:exponent alpha}.
But nevertheless, it is desirable to develop models for cell division processes with alternative division rules as well as life time distributions.
Here, a combination of the results of \cite{schreiber_thae2013_shape_driven, georgii_schreiber_thaele2015} with the present approach could be promising.
Simulation studies for tessellations in a bounded window of the plane $\RR^2$ with several life time rules L and also a variety of division rules D, aimed to adapt  models to real structures, are presented in \cite{leon_ohser_nagel_arscott2019, leon_montero_nagel2023}.
Codes for a simulation of those tessellations 
are available in \cite{crackpatternsimul}.

\bigskip

\section*{Acknowledgment} 
This work was supported by the Center for Mathematical Modeling
ANID Basal FB210005. In particular, Werner Nagel is indebted for the support of his visits at the Center, where in March 2023 this paper could be finished.

\bigskip

\bibliographystyle{APT}
\bibliography{literatur_nagel}

\end{document}